\documentclass{amsart}
\usepackage{amssymb}
\usepackage{amsfonts}

\newtheorem{theorem}{Theorem}[section]
\newtheorem*{theorem*}{Theorem}
\theoremstyle{plain}

\newtheorem{corollary}[theorem]{Corollary}

\newtheorem{definition}[theorem]{Definition}
\newtheorem*{definition*}{Definition}
\newtheorem{example}[theorem]{Example}

\newtheorem{lemma}[theorem]{Lemma}

\newtheorem{proposition}[theorem]{Proposition}

\begin{document}
\title[Infinity-harmonic Maps and Morphisms]{Infinity-harmonic Maps and
Morphisms}
\author{Ye-Lin Ou$^{*}$}
\address{Department of Mathematics\\
Texas A\&M University-Commerce\\
Commerce, Tx. 75429-3011\\
Yelin\_Ou@tamu-commerce.edu}
\author{Tiffany Troutman}
\address{Department of Mathematics\\
Bradley University\\
Peoria, IL 61625\\
ttroutman@bradley.edu}
\author{Frederick Wilhelm}
\address{Department of Mathematics\\
University of California, Riverside\\
Riverside, CA 92521\\
fred@math.ucr.edu}
\thanks{$^*$ Supported by Texas A $\&$ M University-Commerce ``Faculty
Research Enhancement Program" (2008)}
\maketitle

\section{Introduction\protect\medskip}

Distance functions play a major role in our knowledge of Riemannian
geometry. Wherever a distance function is smooth its gradient has constant
norm $1$ and so for trivial reason satisfies the $\infty$-Laplace equation 
\begin{equation*}
\langle \mathrm{grad}\,u,\mathrm{grad}\left\vert \mathrm{grad}u\right\vert
^{2}\rangle =0.
\end{equation*}%
The euclidean version of this equation was introduced by Aronsson \cite{Ar1}
in the 1960's. The solutions are called infinity-harmonic functions.
Geometrically, $u$ is infinity-harmonic if and only if a given integral
curve of its gradient field is parameterized with constant speed.

Since infinity-harmonic functions are a natural generalization of smooth
distance functions, there is great potential for advancing our knowledge of
Riemannian geometry through the lens of this analytic concept. We hope to
lay the ground work for this by answering some basic questions here.
Although there are many interesting examples of infinity-harmonic functions
on Riemannian manifolds, we will work in the broader context of
infinity-harmonic maps.

\begin{definition*}
A $C^{2}$ map $\psi :(M,g)\rightarrow (N,h)$ is said to be infinity-harmonic
if and only if \addtocounter{theorem}{1} 
\begin{equation}
\Delta _{\infty }\left( \psi \right) \equiv \frac{1}{2}d\psi (\mathrm{grad}%
\,|d\psi |^{2})=0,  \label{00E}
\end{equation}%
where 
\begin{equation*}
|d\psi |_{x}^{2}\equiv \sum_{i=1}^{n}h(d\psi (e_{i}),d\psi (e_{i}))
\end{equation*}%
is called the energy density of $\psi $, and $\{e_{i}\}$ is an orthonormal
basis for $T_{x}M$.
\end{definition*}

This generalizes the concept of infinity-harmonic functions on euclidean
space. The definition can also be viewed as the limiting case of the well
known notion of $p$-harmonic map \cite{BG} as $p\rightarrow \infty .$ (See
Proposition \ref{limit as p goes to inf} for details.)

Metric projection (i.e. the closest point map), to an orbit of an isometric
group action is typically not an isometry, or even a Riemannian submersion,
even at the places that it is well defined, smooth, and the orbit is
principle. On the other hand, it is always infinity-harmonic.

\begin{theorem}
\label{isometric grp actions}Let $O\subset M$ be an orbit of an isometric
action by a Lie group $G$ on a Riemannian manifold $M.$ Then metric
projection onto $O$ is an infinity-harmonic map, wherever it is well defined
and smooth.
\end{theorem}

In Section 2, we prove Theorem \ref{isometric grp actions} and give other
examples of infinity-harmonic maps, including projections of multiply warped
products, totally geodesic maps, isometric immersions, Riemannian
submersions, and eigenmaps between spheres.

Section 3 begins with some examples that show that infinity-harmonicity is
not preserved under composition. Motivated by this and the theory of $p$%
-harmonicity, we introduce a subclass of infinity-harmonic maps called
infinity-harmonic morphisms, which preserve solutions to the $\infty$%
-Laplace equation in the following sense.

\begin{definition*}
A map between Riemannian manifolds is said to be an infinity-harmonic
morphism if and only if it pulls back locally defined infinity-harmonic
functions to infinity-harmonic functions.
\end{definition*}

This is motivated by the categorically analogous definition of $p$-harmonic
morphism (\cite{Fu}, \cite{Is}, \cite{Lo}, and \cite{BL}), and is therefore
very appealing. On the other hand, it is a difficult condition to verify.
Fortunately, we will provide an alternative characterization of
infinity-harmonic morphisms that is easier to check. To this end we recall 
\cite{BE}, \cite{BW1}

\begin{definition*}
A map $\varphi :(M,g)\longrightarrow (N,h)$ between Riemannian manifolds is
horizontally weakly conformal with dilation $\lambda :M\longrightarrow
\lbrack 0,\infty )$ if apart from the points where $d\varphi =0$, $d\varphi
_{x}$ is onto and 
\begin{equation*}
h(d\varphi _{x}(X),d\varphi _{x}(Y))=\lambda ^{2}(x)g_{x}(X,Y)
\end{equation*}%
for all horizontal vector fields on $M$.

A horizontally weakly conformal map with dilation $\lambda $ having vertical
gradient is a horizontally homothetic map. A horizontally weakly conformal
map without a critical point is called a horizontally conformal submersion
and a horizontally homothetic map without a critical point is called a
horizontally homothetic submersion.\medskip
\end{definition*}

\begin{theorem}
\label{HWCTHM copy(1)} A map between Riemannian manifolds is an
infinity-harmonic morphism if and only if it is a horizontally weakly
conformal, infinity-harmonic map, and such a map is precisely a horizontally
homothetic map.
\end{theorem}

\medskip

\medskip

This is proven in Section 3. In Section 4, we give several methods to
construct infinity-harmonic maps into Euclidean spaces, characterize those
immersions which are infinity-harmonic maps, and show that isometrically
immersing the target manifold of a map into another manifold does not change
the infinity-harmonicity of the map. Section 5 is devoted to constructions
of infinity-harmonic maps into spheres. We use ideas similar to those of
Smith's in finding harmonic maps into spheres to find infinity-harmonic maps
into spheres by reduction of partial differential equations into ordinary
differential equations. Finally, in Section 6 we examine the effect of a
conformal change on the $\infty$-Laplacian to obtain formulas for the
$\infty$-Laplace equation on spheres and on hyperbolic spaces in terms of
the $\infty$-Laplacian on Euclidean space.

\section{Some Examples and Properties of Infinity-harmonic Maps}

In this section, We will show that metric projection to an orbit of an
isometric group action is always infinity-harmonic, study the relationship
between infinity-harmonic and p-harmonic maps, and give some examples of
infinity-harmonic maps some of which have play important role in
differential geometry.\newline

\subsection{A Metric Projection is Always Infinity-harmonic}

To prove Theorem \ref{isometric grp actions}\textbf{\ }we need the following
Lemma, which can be viewed as a corollary of the Slice Theorem [Theorem 5.4
Bred]. We include a direct proof for the convenience of the reader.

\begin{lemma}
Let $G$ act on $M$ by isometries. Let $O$ be an orbit of $G$, and let $N$ be
a $G$--invariant subset on which the metric projection map $\pi
:N\longrightarrow O$ is defined and smooth. Then the restriction of $\pi $
to any orbit $O_{2}$ contained in $N$ is a submersion (in the smooth sense).
\end{lemma}

\begin{proof}
Let $v\in T_{p}O$ be given, and let $q$ be any point in $\pi ^{-1}\left(
p\right) \cap O_{2}.$ Let $\left\{ g_{t}\right\} _{t>0}\subset G$ be a one
parameter subset so that $\frac{d}{dt}g_{t}\left( p\right) |_{t=0}=v.$ Then
because $\pi $ is the closest point map and the $g_{t}$s are isometries %
\addtocounter{theorem}{1} 
\begin{equation}
\pi \left( g_{t}\left( q\right) \right) =g_{t}\left( p\right) .
\label{symmetries}
\end{equation}%
If $w=\frac{d}{dt}g_{t}\left( q\right) |_{t=0},$ it follows that 
\begin{eqnarray*}
d\pi \left( w\right) &=&\frac{d}{dt}\pi \left( g_{t}\left( q\right) \right)
|_{t=0} \\
&=&v,
\end{eqnarray*}%
so $\pi |_{O_{2}}$ is a submersion.
\end{proof}

The same proof also gives us part i of the following lemma.

\begin{lemma}
\label{local submersion}Let $G$ act on $M$ by isometries. Let $O$ be an
orbit of $G$, and let $N$ be a open subset on which the metric projection
map $\pi :N\longrightarrow O$ is defined and smooth. For $q\in N$ let $%
\mathcal{O}_{q}$ be the orbit through $q.$

\begin{description}
\item[i] For any $q\in N,$ the differential of $\pi |_{\mathcal{O}_{q}}$ is
onto, and

\item[ii] $T_{q}\mathcal{O}_{q}$ contains the horizontal distribution, $%
\mathcal{H}_{q}$, of $\pi $ at $q.$
\end{description}
\end{lemma}

\noindent \textbf{Proof of Part ii: }Given $q\in N,$ let $\mathcal{V}_{q}$
be the kernel of $d\pi _{q},$ and let $T\mathcal{O}^{\perp }$ be the
orthogonal complement of $\mathcal{V}_{q}\cap T\mathcal{O}$ in $T\mathcal{O}%
. $ It follows from dimension counting and part $i$ of the lemma that $T%
\mathcal{O}^{\perp }$ coincides with $\mathcal{H}_{q},$ the orthogonal
complement of $\mathcal{V}_{q}$ in $T_{q}M.$ $\square $

\begin{corollary}
\label{G-familly}Let $G$ act on $M$ by isometries . Let $\mathcal{O}$ be an
orbit of $G$, and let $N$ be a open subset on which the metric projection
map $\pi :N\longrightarrow \mathcal{O}$ is defined and smooth. Then for any
smooth $\pi $--horizontal curve, $\gamma ,$ in $N$ that there is a smooth
one parameter subset of isometries $\left\{ g_{t}\right\} \subset G$ so that%
\begin{equation*}
g_{t}\gamma \left( 0\right) =\gamma \left( t\right) .
\end{equation*}
\end{corollary}

\textbf{Proof of Theorem \ref{isometric grp actions}: }Let metric
projection, $\pi :N\longrightarrow O$ be defined and smooth on the subset $N$
of $M.$

Let $\gamma $ be a horizontal curve for $\pi .$ From Corollary \ref%
{G-familly}, there is a one parameter subset of isometries $\left\{
g_{t}\right\} \subset G$ so that

\begin{equation*}
g_{t}\gamma \left( 0\right) =\gamma \left( t\right) .
\end{equation*}
It follows from equation (\ref{symmetries}) that $G$ acts by symmetries of $%
\pi .$ In particular, it preserves the horizontal and vertical distributions
of $\pi ,$ it follows that 
\begin{equation*}
\left\vert d\pi \right\vert ^{2}\left( x\right) =\left\vert d\pi \right\vert
^{2}\left( g_{t}x\right)
\end{equation*}%
for all $t.$ In particular, $\mathrm{grad}|d\pi |^{2}$ is vertical. $\square$%
\newline

\subsection{Relationship Between Infinity-harmonic and $p$-harmonic Maps}

Recall that, a \emph{$p$-harmonic map} ($p>1$) is a map $\varphi
:(M,g)\longrightarrow (N,h)$ between Riemannian manifolds such that $\varphi
|\Omega $ is a critical point of the $p$-energy 
\begin{equation}
E_{p}\left( \varphi ,\Omega \right) =\frac{1}{p}{\int }_{\Omega }\left\vert 
\mathrm{d}\varphi \right\vert ^{p}\mathrm{d}x  \notag
\end{equation}%
for every compact subset $\Omega $ of $M$. Locally, $p$-harmonic maps are
solutions of the following systems of PDEs: \addtocounter{theorem}{1} 
\begin{equation}
\Delta _{p}(\varphi )={\left\vert \mathrm{d}\varphi \right\vert }%
^{p-2}\Delta _{2}(\varphi )+(p-2){\left\vert \mathrm{d}\varphi \right\vert }%
^{p-4}\mathrm{\ d}\varphi (\mathrm{grad}{\left\vert \mathrm{d}\varphi
\right\vert })=0,  \label{phm}
\end{equation}%
where $\Delta _{2}(\varphi )=\mathrm{Trace}_{g}\nabla \mathrm{d}\varphi $
denotes the tension field of $\varphi $. When $p=2$, we get the familiar
notion of harmonic maps which include geodesics, harmonic functions, minimal
isometric immersions, and Riemannian submersions with minimal fibers as
special cases (See \cite{EL1}, \cite{EL2}, \cite{EL3}, and \cite{SY}). 
\newline

We point out that the definition of infinity-harmonic map can be viewed as
the limiting case of the notion of $p$--harmonic map as $p\rightarrow \infty 
$ in the following sense.

\begin{proposition}
\label{limit as p goes to inf}For any $C>0,$ 
\begin{equation*}
\lim_{p\rightarrow \infty }\sup_{\varphi \in H_{C}^{p}}\Delta _{\infty
}\left( \varphi \right) =0
\end{equation*}%
where $H_{C}^{p}$ is the class of all $p$--harmonic maps $\varphi $ with $%
\left\vert \mathrm{d}\varphi \right\vert ^{2}\left\vert \Delta _{2}(\varphi
)\right\vert \leq C.$
\end{proposition}

\begin{proof}
Dividing the $p$-harmonic equation by $(p-2){\left\vert \mathrm{d}\varphi
\right\vert }^{p-4}$ gives \addtocounter{theorem}{1} 
\begin{equation}
\frac{\left\vert \mathrm{d}\varphi \right\vert ^{2}\,\Delta _{2}(\varphi )}{%
(p-2)}+\frac{1}{2}\mathrm{d}\varphi (\mathrm{grad}\left\vert \mathrm{d}%
\varphi \right\vert ^{2})=0.
\end{equation}%
So within the class of $p$--harmonic maps $\varphi $ with $\left\vert 
\mathrm{d}\varphi \right\vert ^{2}\left\vert \Delta _{2}(\varphi
)\right\vert \leq C,$ we can make the $\left\vert \frac{1}{2}\mathrm{d}%
\varphi (\mathrm{grad}\left\vert \mathrm{d}\varphi \right\vert
^{2})\right\vert $ as small as we please be letting $p\rightarrow \infty .$
\end{proof}

Another relationship between $p$-harmonic and infinity-harmonic maps is a
consequence of equation \ref{phm}.

\begin{proposition}
\label{RELA} If a map is $p$-harmonic for two different $p$ values, then it
is infinity-harmonic; An infinity-harmonic map is also a harmonic map if
and only if it is a $p$-harmonic map for any $p\neq 2$.
\end{proposition}

\subsection{Some Examples of Infinity-harmonic Maps}

Besides metric projections, the following important and familiar classes of
maps are infinity-harmonic.

\begin{example}
{$[\mathbf{Infinity-harmonic}$ $\mathbf{{functions}]}$} A real-valued
function 
\begin{equation*}
u:(M,g)\longrightarrow \mathbb{R}
\end{equation*}%
on a Riemannian manifold is infinity-harmonic if and only if $u$ is a
solution of $\infty$-Laplace equation : \addtocounter{theorem}{1} 
\begin{eqnarray}
\Delta _{\infty }u &=&\frac{1}{2}\mathrm{d}u(\mathrm{grad}\,\left\vert 
\mathrm{\ grad}\,u\right\vert ^{2})  \label{00f} \\
&=&\frac{1}{2}g(\mathrm{grad}\,u,\mathrm{grad}\,\left\vert \mathrm{grad}%
\,u\right\vert ^{2})  \notag \\
&=&0.  \notag
\end{eqnarray}

For $u:\Omega \subset \mathbb{R}^{m}\longrightarrow \mathbb{R}$, this
becomes Aronsson's $\infty$-Laplace equation.
\end{example}

\begin{example}
{$[\mathbf{Maps\,with\,constant\,energy\,density}]$} Any map with constant
energy density, 
\begin{equation*}
\left\vert \mathrm{d}\varphi \right\vert ^{2}=\mathrm{constant},
\end{equation*}
is infinity-harmonic. This class includes

\begin{itemize}
\item any totally geodesic map between Riemannian manifolds. Recall that a
map $\varphi :(M^{m},g)\longrightarrow (N^{n},h)$ is totally geodesic if its
second fundamental form vanishes identically, i.e., $\nabla \mathrm{d}%
\varphi =0$. It is not difficult to see that $\varphi $ is totally geodesic
if and only if it carries geodesics to geodesics. It is well known \cite{ER}
that a totally geodesic map has constant rank and constant energy density;

\item any eigenmap between spheres $\varphi :S^{m}\longrightarrow S^{n}$.
Recall that an eigenmap is a harmonic map between spheres with constant
energy density, which can be characterized as the restriction to $S^{m}$ of
a map $F:\mathbb{R}^{m+1}\longrightarrow \mathbb{R}^{n+1}$ whose components
are harmonic homogeneous polynomial of a common degree \cite{ER};

\item the globally defined nonlinear complex-valued functions $\varphi :%
\mathbb{R}^{m}\longrightarrow \mathbb{R}^{2}\equiv \mathbb{C}$ with $\varphi
(x_{1},\ldots ,x_{m})=\lambda _{1}e^{ix_{1}}+\ldots +\lambda _{m}e^{ix_{m}}$
, where $\lambda _{k}\,(k=1,\ldots ,m)$ are constant real numbers. One can
easily check that this map has constant energy density $\left\vert \mathrm{d}%
\varphi \right\vert ^{2}=\sum_{k=1}^{m}|\lambda _{k}|^{2}$. Note that a map
of this class does not belong to any of the above classes, for instance, $%
\varphi :\mathbb{R}^{3}\longrightarrow \mathbb{R}^{2}$ given by $\varphi
(x,y,z)=(\cos x+\cos y+\cos z,\;\sin x+\sin y+\sin z)$ is a globally defined
smooth nonlinear infinity-harmonic map which is neither an isometric
immersion nor a Riemannian submersion;

\item Let $\varphi :\mathbb{R}^{m}\longrightarrow \mathbb{R}^{m}$ be defined
by 
\begin{eqnarray}
&&\hskip3cm\varphi (x_{1},\ldots ,x_{m})=  \notag \\
&&(\cos x_{1}+\sin x_{2},\;\cos x_{2}+\sin x_{3},\;\ldots ,\;\cos
x_{m-1}+\sin x_{m},\;\cos x_{m}+\sin x_{1}).  \notag
\end{eqnarray}%
A straightforward computation gives the energy density $\left\vert \mathrm{d}%
\varphi \right\vert ^{2}=m$.
\end{itemize}
\end{example}

\begin{example}
{$[\mathbf{Infinity-harmonic\,curves}]$} Any regular curve $\gamma
:(a,b)\longrightarrow (M^{m},g)$ is an infinity-harmonic map provided it is
parametrized by arc length.
\end{example}

The following example provides a large class of infinity-harmonic maps with
nonconstant energy density.

\begin{example}
\label{Multi-warped products}{$[\mathbf{Projection\,of\,multiply\,warped%
\,products}]$} Recall that a multiply warped product of Riemannian manifolds 
$(B,g_{B})$ and $(F_{1},h_{1}),\ldots ,(F_{k},h_{k})$ is the smooth manifold 
$M=B\times F_{1}\times \ldots \times F_{k}$ with the metric 
\begin{equation*}
g_{B}+\lambda _{1}^{2}h_{1}+\ldots +\lambda _{k}^{2}h_{k},
\end{equation*}%
where $\lambda _{1},\ldots ,\lambda _{k}:B\longrightarrow (0,\infty )$ are
called warping functions. We denote the resulting Riemannian manifolds by $%
B\times _{\lambda _{1}^{2}}F_{1}\times \ldots \times _{\lambda
_{k}^{2}}F_{k} $.

Let 
\begin{equation*}
\pi :B\times _{\lambda _{1}^{2}}F_{1}\times \ldots \times _{\lambda
_{k}^{2}}F_{k}\longrightarrow (F_{1}\times \ldots \times F_{k},h_{1}+\ldots
+h_{k})
\end{equation*}%
be projection. A simple computation gives that the energy density of $\pi $
is 
\begin{equation*}
\left\vert \mathrm{d}\pi \right\vert ^{2}=\lambda _{1}^{-2}+\ldots +\lambda
_{k}^{-2}.
\end{equation*}
Since the gradients of all of the $\lambda _{i}$s are tangent to the
\textquotedblleft B--factors\textquotedblright\ they are all vertical for $%
\pi ,$ and $\pi $ is infinity-harmonic. In particular we have

\begin{itemize}
\item the projection $\pi :(\mathbb{R}^{3},g_{Sol})\longrightarrow (\mathbb{R%
}^{2},\mathrm{d}x^{2}+\mathrm{d}y^{2})$ with $\pi (x,y,z)=(x,y)$ is an
infinity-harmonic map, where $(\mathbb{R}^{3},g_{Sol})$ denotes the
Sol space, one of Thurston's eight $3$-dimensional geometries, which can be
viewed as 
\begin{equation*}
(\mathbb{R}^{3},g_{Sol})=(\mathbb{R}\times \mathbb{R}\times \mathbb{R},e^{2z}%
\mathrm{d}x^{2}+e^{-2z}\mathrm{d}y^{2}+\mathrm{d}z^{2}).
\end{equation*}

\item the projection from $3$-sphere onto the Clifford torus 
\begin{equation*}
\varphi :S^{3}\setminus \{\Gamma _{1},\Gamma _{2}\}\equiv ((0,\frac{\pi }{2}%
)\times S^{1}\times S^{1},\;\mathrm{dt}^{2}+\sin ^{2}t\,\mathrm{d\theta _{1}}%
^{2}+\cos ^{2}t\,\mathrm{d\theta _{2}}^{2})\longrightarrow S^{1}\times S^{1}
\end{equation*}%
with $\varphi (t,\theta _{1},\theta _{2})=(\theta _{1},\theta _{2})$ is an
infinity-harmonic submersion with nonconstant energy density $|\mathrm{d}%
\varphi |^{2}=\frac{1}{\sin ^{2}t}+\frac{1}{\cos ^{2}t}$.
\end{itemize}
\end{example}

\begin{proposition}
\label{P1} \label{inf harm subm}A submersion $\pi :(M^{m},g)\longrightarrow
(N^{n},h)$ is infinity-harmonic if and only if the gradient of the energy
density of $\pi $ is vertical. In particular, a horizontally conformal
submersion is infinity-harmonic if and only if it is a horizontally
homothetic submersion.
\end{proposition}

\section{Infinity-harmonic Morphisms}

As the following example shows, infinity-harmonicity is not preserved under
composition of infinity-harmonic maps.

\begin{example}
The linear map $\phi :\mathbb{R}^{2}\setminus \{0\}\longrightarrow \mathbb{R}%
^{2}\setminus \{0\}$ with $\phi (x,y)=(x,2y)$ is an infinity-harmonic map
since it has constant energy density. It is well known that $f:\mathbb{R}%
^{2}\setminus \{0\}\longrightarrow \mathbb{R}$ with $f(x,y)=\sqrt{%
x^{2}+y^{2}\;}$ is an infinity-harmonic function. One can easily check that
the pull-back function $\phi ^{\ast }f=\sqrt{x^{2}+4y^{2}\;}$ is not an
infinity-harmonic function.\medskip
\end{example}

In this section we will prove Theorem \ref{HWCTHM copy(1)}. Using it we see
that the following are examples of infinity-harmonic morphisms.

\begin{example}
Riemannian Submersions are infinity-harmonic morphisms.
\end{example}

\begin{example}
The projection of a warped product onto the base is a Riemannian submersion.
The projection onto the fiber is a horizontally homothetic submersion and
hence an infinity-harmonic morphism.
\end{example}

\begin{example}
Radial projection of $R^{n+1}\setminus \{0\}$ onto $S^{n}$ is an infinity-harmonic morphism. In this example, $\lambda (x)=\frac{1}{|x|}$.
\end{example}

\begin{example}
Infinity-harmonic functions on Riemannian manifolds can be viewed as
horizontally weakly conformal, infinity-harmonic maps, and hence are
infinity-harmonic morphisms.
\end{example}

We refer the reader to \cite{OW} for other examples of horizontally
homothetic submersions.

We prove Theorem \ref{HWCTHM copy(1)} with three lemmas. The first of which
is as follows.

\begin{lemma}
\label{toHWC}If $\pi :E\rightarrow B$ is an infinity-harmonic morphism, then 
$\pi $ is a horizontally weakly conformal map.
\end{lemma}

Before proving this we study the linear case.

\begin{proposition}
\label{linear inf harm morph} A linear map $\varphi :\mathbb{R}%
^{n+k}\longrightarrow \mathbb{R}^{n}$ is an infinity-harmonic morphism if
and only if it is a horizontally conformal surjective submersion. In other
words, it can be written as the composition of a homothety, an isometry, and
an orthogonal projection.

In fact, if $\varphi $ is onto but not horizontally weakly conformal, then $%
\varphi ^{\ast }\left( \mathrm{dist}\left( 0,\cdot \right) \right) $ is not
infinity-harmonic and 
\begin{equation*}
\lim \sup_{p\rightarrow 0}\left\vert \Delta _{\infty }\left[ \varphi ^{\ast
}\left( \mathrm{dist}\left( 0,\cdot \right) \right) \right] |_{p}\right\vert
=\infty .
\end{equation*}
\end{proposition}

\begin{proof}
First we consider the case when $\varphi :\mathbb{R}^{n+k}\longrightarrow 
\mathbb{R}^{n}$ is onto. We have%
\begin{eqnarray*}
\left\langle \nabla \left\vert \nabla \left( f\circ \varphi \right)
\right\vert ^{2},e_{\alpha }\right\rangle &=&d\left\vert \nabla \left(
f\circ \varphi \right) \right\vert ^{2}\left[ e_{\alpha }\right] \\
&=&\displaystyle\sum_{i}d\left( \left\langle \nabla f,d\varphi \left(
e_{i}\right) \right\rangle ^{2}\right) \left[ e_{\alpha }\right] \\
&=&2\displaystyle\sum_{i}\left\langle \nabla f,d\varphi \left( e_{i}\right)
\right\rangle \left[ d\left\langle \nabla f,d\varphi \left( e_{i}\right)
\right\rangle \right] \left[ e_{\alpha }\right] \\
&=&2\displaystyle\sum_{i}\left\langle \nabla f,d\varphi \left( e_{i}\right)
\right\rangle \left[ \left\langle \nabla _{e_{\alpha }}\left( \nabla f\circ
\varphi \right) ,d\varphi \left( e_{i}\right) \right\rangle +\left\langle
\nabla f,\nabla _{e_{\alpha }}d\varphi \left( e_{i}\right) \right\rangle %
\right] ,
\end{eqnarray*}

\noindent where $\left( \nabla f\circ \varphi \right) ,\nabla _{e_{\alpha
}}\nabla f,d\varphi \left( e_{i}\right) ,$and $\nabla _{e_{\alpha }}d\varphi
\left( e_{i}\right) $ are being viewed as vector fields along $\varphi .$

Since $\varphi $ is linear, $d\varphi \left( e_{i}\right) $ is a constant
vector field, and the second term vanishes. So 
\begin{eqnarray*}
\left\langle \nabla \left\vert \nabla \left( f\circ \varphi \right)
\right\vert ^{2},e_{\alpha }\right\rangle &=&2\displaystyle%
\sum_{i}\left\langle \nabla f,d\varphi \left( e_{i}\right) \right\rangle %
\left[ \left\langle \nabla _{e_{\alpha }}\left( \nabla f\circ \varphi
\right) ,d\varphi \left( e_{i}\right) \right\rangle \right] \text{ so} \\
\nabla \left\vert \nabla \left( f\circ \varphi \right) \right\vert ^{2} &=&2%
\displaystyle\sum_{\alpha ,i}\left\langle \nabla f,d\varphi \left(
e_{i}\right) \right\rangle \left[ \left\langle \nabla _{e_{\alpha }}\left(
\nabla f\circ \varphi \right) ,d\varphi \left( e_{i}\right) \right\rangle %
\right] e_{\alpha }
\end{eqnarray*}%
Combining this with 
\begin{equation*}
\nabla \left( f\circ \varphi \right) =\displaystyle\sum_{\alpha
}\left\langle \nabla f,d\varphi \left( e_{a}\right) \right\rangle e_{a}
\end{equation*}%
we get \addtocounter{theorem}{1} 
\begin{eqnarray}
\Delta _{\infty }\left[ \varphi ^{\ast }\left( f\right) \right] &=&\frac{1}{2%
}\left\langle \nabla \left\vert \nabla \left( f\circ \varphi \right)
\right\vert ^{2},\nabla \left( f\circ \varphi \right) \right\rangle  \notag
\label{inf-lapl of pull back} \\
&=&\displaystyle\sum_{\alpha ,i}\left\langle \nabla f,d\varphi \left(
e_{i}\right) \right\rangle \left[ \left\langle \nabla _{e_{\alpha }}\left(
\nabla f\circ \varphi \right) ,d\varphi \left( e_{i}\right) \right\rangle %
\right] \left\langle \nabla f,d\varphi \left( e_{a}\right) \right\rangle
\label{inf
lapl of Pull Back}
\end{eqnarray}

Now set 
\begin{equation*}
f\left( x\right) =\sqrt{x_{1}^{2}+x_{2}^{2}+\cdots +x_{n}^{2}}
\end{equation*}%
Then 
\begin{equation*}
\nabla f=\frac{1}{f}\left( x_{1},x_{2},\ldots ,x_{n}\right)
\end{equation*}

Since $f$ is a distance function we have 
\begin{equation*}
\nabla _{\nabla f}\nabla f\equiv 0.
\end{equation*}%
Since $f$ is the distance function from the origin, we also have 
\begin{equation*}
\nabla _{z}\nabla f=\frac{z}{f}
\end{equation*}%
for all $z\perp \nabla f.$

To evaluate $\Delta _{\infty }\left( f\circ \varphi \right) $ at $\tilde{p}%
\in \mathbb{R}^{n+k}$ using \ref{inf lapl of Pull Back}, choose an
orthonormal basis $\left\{ v_{0},v_{1},v_{2},\ldots ,v_{n+k-1}\right\} $ for 
$T_{\tilde{p}}\mathbb{R}^{n+k}$ so that 
\begin{equation*}
\left\langle d\varphi \left( v_{i}\right) ,\nabla f\right\rangle |_{\varphi
\left( \tilde{p}\right) }=0,\text{ for }i=1,2,3,\ldots ,n+k-1
\end{equation*}

Then using \ref{inf lapl of Pull Back} 
\begin{eqnarray*}
\Delta _{\infty }\left( f\circ \varphi \right) _{\tilde{p}} &=&\displaystyle%
\sum_{\alpha ,i}\left\langle \nabla f,d\varphi \left( v_{i}\right)
\right\rangle \left[ \left\langle \nabla _{v_{\alpha }}\left( \nabla f\circ
\varphi \right) ,d\varphi \left( v_{i}\right) \right\rangle \right]
\left\langle \nabla f,d\varphi \left( v_{a}\right) \right\rangle \\
&=&\left\langle \nabla f,d\varphi \left( v_{0}\right) \right\rangle ^{2} 
\left[ \left\langle \nabla _{v_{0}}\left( \nabla f\circ \varphi \right)
,d\varphi \left( v_{0}\right) \right\rangle \right]
\end{eqnarray*}%
Since \textrm{span}$\left\{ v_{1},\ldots ,v_{n+k-1}\right\} $ contains the
vertical space for $\varphi ,$ $v_{0}$ is horizontal for $\varphi .$ It
follows that 
\begin{equation*}
\Delta _{\infty }\left( f\circ \varphi \right) _{\tilde{p}}=\left\langle
\nabla f,d\varphi \left( v_{0}\right) \right\rangle ^{2}\left[ \left\langle
\nabla _{d\varphi \left( v_{0}\right) }\nabla f,d\varphi \left( v_{0}\right)
\right\rangle \right] .
\end{equation*}%
Since $\nabla _{d\varphi \left( v_{0}\right) }\nabla f$ is proportional to
the component, $d\varphi \left( v_{0}\right) ^{\perp },$ of $d\varphi \left(
v_{0}\right) $ that is perpendicular to $\nabla f,$ it follows that $\Delta
_{\infty }\left( f\circ \varphi \right) _{\tilde{p}}=0$ if and only if $%
d\varphi \left( v_{0}\right) $ is proportional to $\nabla f|_{p}.$ This is
equivalent to saying that $d\varphi $ maps the orthogonal spaces $\mathrm{%
span}\left\{ v_{0}\right\} $ and $\mathrm{span}\left\{ v_{1},v_{2},\ldots
,v_{n+k-1}\right\} $ to the orthogonal spaces $\mathrm{span}\left\{ \nabla
f|_{\varphi \left( \tilde{p}\right) }\right\} $ and $\mathrm{span}\left\{
\nabla f|_{\varphi \left( \tilde{p}\right) }\right\} ^{\perp }.$ By varying $%
p\equiv \varphi \left( \tilde{p}\right) $ we can make $\nabla f|_{p}$ point
in any direction, and it follows that $d\varphi =\varphi $ preserves all
angles in its horizontal space. So an onto linear infinity-harmonic morphism
is a weakly conformal (surjective) submersion as claimed.

Now suppose $\varphi :\mathbb{R}^{n+k}\longrightarrow \mathbb{R}^{n}$ is a
linear infinity-harmonic morphism, that is not onto. We may post compose
with an orthogonal transformation $\mathbb{R}^{n}\longrightarrow \mathbb{R}%
^{n}$ to obtain a a linear infinity-harmonic morphism whose image is
contained in the coordinate subspace $\mathbb{R}^{m}\times \left\{ 0\right\}
\subset \mathbb{R}^{m}\times \mathbb{R}^{n-m}.$ Applying the result just
proven to $\varphi :\mathbb{R}^{n+k}\longrightarrow \mathbb{R}^{m}$ we see
that $\varphi $ is a horizontally conformal linear submersion.

Let $f:\mathbb{R}^{n}\longrightarrow \mathbb{R}^{n}$ be the distance
function from $\left( 0,\ldots ,0,\varepsilon \right) \notin \mathrm{Im}%
\left( \varphi \right) .$ Then the curve 
\begin{equation*}
c:t\longmapsto \left( t,0,\ldots 0\right)
\end{equation*}%
has the same image as an integral curve, $\gamma$, of $\nabla f|_{\mathbb{R}%
^{k}\times \left\{ 0\right\} }$, only the velocity field of $\gamma $ at $%
c\left( t\right) $ is%
\begin{equation*}
\left( \frac{t}{\sqrt{t^{2}+\varepsilon }},0,0,\ldots ,0\right)
\end{equation*}

Notice in particular that this integral curve of $\nabla f|_{\mathbb{R}%
^{k}\times \left\{ 0\right\} }$ with variable speed, and that $\nabla \left(
\varphi ^{\ast }\left( f\right) \right) $ is a horizontal lift of $\nabla
f|_{\mathbb{R}^{k}\times \left\{ 0\right\} }$. Since $\varphi :\mathbb{R}%
^{n+k}\longrightarrow $ $\mathbb{R}^{k}\times \left\{ 0\right\} \subset 
\mathbb{R}^{k}\times \mathbb{R}^{n-k}$ is a horizontally conformal linear
submersion it follows that $\nabla \left( \varphi ^{\ast }\left( f\right)
\right) $ will also have an integral curve that is parameterized with
variable speed. So $\varphi ^{\ast }\left( f\right) $ is not infinity-harmonic and $\varphi $ is not an infinity-harmonic morphism.

To estimate the $\infty$-Laplacian of $\varphi ^{\ast }\left( f\right) $ in
this event, we compute%
\begin{eqnarray*}
\left\vert \nabla f|_{\mathbb{R}^{k}\times \left\{ 0\right\} }\right\vert
^{2} &=&\frac{t^{2}}{t^{2}+\varepsilon } \\
\nabla \left\vert \nabla f|_{\mathbb{R}^{k}\times \left\{ 0\right\}
}\right\vert ^{2} &=&\left( \frac{2t\left( t^{2}+\varepsilon \right) -2tt^{2}%
}{\left( t^{2}+\varepsilon \right) ^{2}},0,\ldots ,0\right) \\
&=&\left( \frac{2t\varepsilon }{\left( t^{2}+\varepsilon \right) ^{2}}%
,0,\ldots ,0\right) \\
\left\langle \nabla \left\vert \nabla f|_{\mathbb{R}^{k}\times \left\{
0\right\} }\right\vert ^{2},\nabla f|_{\mathbb{R}^{k}\times \left\{
0\right\} }\right\rangle &=&\frac{2t^{2}\varepsilon }{\left(
t^{2}+\varepsilon \right) ^{5/2}}
\end{eqnarray*}%
So if $t^{2}=\varepsilon $ we have 
\begin{eqnarray*}
\left\langle \nabla \left\vert \nabla f|_{\mathbb{R}^{k}\times \left\{
0\right\} }\right\vert ^{2},\nabla f|_{\mathbb{R}^{k}\times \left\{
0\right\} }\right\rangle &=&O\left( \frac{\varepsilon ^{2}}{\varepsilon
^{5/2}}\right) \\
&=&O\left( \frac{1}{\varepsilon ^{1/2}}\right)
\end{eqnarray*}

So if $\varphi $ is linear, horizontally weakly conformal, and not onto,
then we can find infinity-harmonic functions on the target that pull back to
functions with arbitrarily large $\infty$-laplacians.

Conversely, if $\varphi :\mathbb{R}^{n+k}\longrightarrow \mathbb{R}^{n}$ is
an linear horizontally conformal surjective submersion and $f:\mathbb{R}%
^{n}\longrightarrow \mathbb{R}$ is infinity-harmonic, then $\nabla \left(
\varphi ^{\ast }\left( f\right) \right) $ is obtained as a horizontal lift
of $\nabla f.$ Since the integral curves of $\nabla f$ are parameterized
with constant speed and $\varphi :\mathbb{R}^{n+k}\longrightarrow \mathbb{R}%
^{n}$ is a linear horizontally conformal submersion, it follows that the
integral curves of $\nabla \left( \varphi ^{\ast }\left( f\right) \right) $
are also parameterized with constant speed and hence that $\varphi ^{\ast
}\left( f\right) $ is infinity-harmonic and that $\varphi $ is an infinity-harmonic morphism.

If $\varphi $ is onto but not horizontally weakly conformal, then as we have
seen $\Delta _{\infty }\left( f\circ \varphi \right) _{\tilde{p}}\neq 0$ for
some $\tilde{p},$ and 
\begin{eqnarray*}
\Delta _{\infty }\left( f\circ \varphi \right) _{\tilde{p}} &=&\left\langle
\nabla f,d\varphi \left( v_{0}\right) \right\rangle ^{2}\left[ \left\langle
\nabla _{d\varphi \left( v_{0}\right) }\nabla f,d\varphi \left( v_{0}\right)
\right\rangle \right] \\
&=&\left\langle \nabla f,d\varphi \left( v_{0}\right) \right\rangle ^{2} 
\left[ \left\langle \nabla _{d\varphi \left( v_{0}\right) ^{\perp }}\nabla
f,d\varphi \left( v_{0}\right) ^{\perp }\right\rangle \right]
\end{eqnarray*}

where we have again used the fact that $\nabla _{d\varphi \left(
v_{0}\right) }\nabla f$ is proportional to $d\varphi \left( v_{0}\right)
^{\perp }.$

Therefore%
\begin{equation*}
\Delta _{\infty }\left( f\circ \varphi \right) _{\tilde{p}}=\frac{1}{f}%
\left\langle \nabla f,d\varphi \left( v_{0}\right) \right\rangle
^{2}\left\vert d\varphi \left( v_{0}\right) \right\vert ^{2}\sin
^{2}\sphericalangle \left( d\varphi \left( v_{0}\right) ,\nabla f\right)
\end{equation*}

Letting $\tilde{p}$ approach the origin along a radial line, all quantities
on the right hand side stay fixed, except, $\frac{1}{f}$ which goes to $%
\infty .$ It follows that 
\begin{equation*}
\lim \sup_{\tilde{p}\rightarrow 0}\left\vert \Delta _{\infty }\left( f\circ
\varphi \right) _{\tilde{p}}\right\vert =\infty .
\end{equation*}
\end{proof}

To prove Lemma \ref{toHWC} we combine our characterization of linear
infinity-harmonic morphisms with the principle that the set of infinity-harmonic functions is closed in the $C^{1}$--topology.

This principle is embodied in the next two propositions.

\begin{proposition}
\label{C^1 close} The set of infinity-harmonic functions is closed in the $%
C^{1}$--topology. I.e. if $f:U\longrightarrow \mathbb{R}$ is any locally
defined real valued $C^{2}$--function on $M$ that is not infinity-harmonic,
then there is an $\epsilon >0$ so that if $h$ is any $C^{2}$--function with%
\begin{eqnarray*}
\left\vert f\left( x\right) -h\left( x\right) \right\vert &<&\epsilon ,\text{
and} \\
\left\vert df\left( v\right) -dh\left( v\right) \right\vert &<&\epsilon
\end{eqnarray*}%
for all $x\in U$ and all unit vectors $v\in TM$, then $h$ is not infinity-harmonic.
\end{proposition}

\begin{proof}
Since $f$ is not infinity-harmonic, there is an integral curve $\gamma :%
\left[ a,b\right] \longrightarrow U$ of $\nabla f$ and a constant $M>0$ so
that \addtocounter{theorem}{1} 
\begin{equation}
\left\vert \left\vert \nabla f_{\gamma \left( b\right) }\right\vert
-\left\vert \nabla f_{\gamma \left( a\right) }\right\vert \right\vert
>M\left( b-a\right) .
\end{equation}

By continuity, there is a neighborhood $V_{a}$ of $\gamma \left( a\right) $
so that any integral curve $c$ of $\nabla f$ passing through $V_{a}$ can be
parameterized on $\left[ a,b\right] $ and satisfies 
\begin{eqnarray*}
\left\vert \left\vert \nabla f_{c\left( b\right) }\right\vert -\left\vert
\nabla f_{c\left( a\right) }\right\vert \right\vert &>&\frac{M}{2}\left(
b-a\right) \\
\left\vert \left\vert \nabla f_{c\left( t\right) }\right\vert -\left\vert
\nabla f_{\gamma \left( t\right) }\right\vert \right\vert &\leq &\frac{M}{100%
}\left( b-a\right) .
\end{eqnarray*}%
Let $\Phi $ be the flow of $\nabla f$ and set 
\begin{equation*}
V=\cup _{t\in \left[ a,b\right] }\Phi _{t}\left( V_{a}\right) .
\end{equation*}%
Set $\epsilon <\frac{M}{100}\left( b-a\right) ,$ and choose $h$ so that 
\begin{eqnarray*}
\left\vert f\left( x\right) -h\left( x\right) \right\vert &<&\epsilon ,\text{
and} \\
\left\vert df\left( v\right) -dh\left( v\right) \right\vert &<&\epsilon
\end{eqnarray*}%
for all $x\in U$ and all unit vectors $v\in TM.$ Require further that $h$ is
close enough to $f$ in the $C^{1}$--topology so that the integral curve $%
\beta $ of $\nabla h$ that starts at $\gamma \left( a\right) $ is
parameterized on $\left[ a,b\right] $ and stays in $V.$ Then 
\begin{eqnarray*}
\left\vert \left\vert \nabla h_{\beta \left( b\right) }\right\vert
-\left\vert \nabla h_{\beta \left( a\right) }\right\vert \right\vert &\geq
&\left\vert \left\vert \nabla f_{\beta \left( b\right) }\right\vert
-\left\vert \nabla f_{\beta \left( a\right) }\right\vert \right\vert
-2\epsilon \\
&\geq &\left\vert \left\vert \nabla f_{\gamma \left( b\right) }\right\vert
-\left\vert \nabla f_{\gamma \left( a\right) }\right\vert \right\vert
-4\epsilon \\
&\geq &\frac{M}{4}\left( b-a\right) .
\end{eqnarray*}%
So $h$ is not infinity-harmonic.
\end{proof}

\begin{definition}
We call a map $\Psi :M\longrightarrow N$ and $\epsilon $-isometry provided, 
\begin{equation*}
\left\vert \left\vert d\Psi \left( v\right) \right\vert -1\right\vert
<\epsilon
\end{equation*}%
for all unit vectors $v.$
\end{definition}

An argument similar to the previous proposition gives us

\begin{proposition}
Let $f$ be a locally defined real valued $C^{2}$--function on $M$ that is
not infinity-harmonic. There is an $\epsilon >0$ so that if $\Psi $ is an $%
\epsilon $--isometry$,$ then $\Psi ^{\ast }\left( f\right) $ is not infinity-harmonic.

More generally, if $\pi :E\longrightarrow B$ is a $C^{2}$--map and $%
f:B\longrightarrow \mathbb{R}$ is such that $\pi ^{\ast }\left( f\right) $
is not infinity-harmonic, then there is an $\epsilon >0$ so that if $\Phi :$ 
$E\longrightarrow E$ and $\Psi :B\longrightarrow B$ are $\epsilon $%
--isometries$,$ then $\left( \Psi \circ \pi \circ \Phi \right) ^{\ast
}\left( f\right) $ is not infinity-harmonic.
\end{proposition}

We can now offer the proof of Lemma \ref{toHWC}.

\begin{proof}
Suppose $\pi :E\longrightarrow B$ is a map that is not horizontal weakly
conformal. Let $\pi \left( \tilde{p}\right) =p$ and suppose that $d\pi _{%
\tilde{p}}$ is onto, but not horizontal weakly conformal. Let 
\begin{eqnarray*}
f &:&T_{\tilde{p}}E\longrightarrow \mathbb{R}\text{ be} \\
f &=&\mathrm{dist}\left( 0,\cdot \right) \circ d\pi _{\tilde{p}}
\end{eqnarray*}%
and 
\begin{eqnarray*}
h &:&T_{\tilde{p}}E\longrightarrow \mathbb{R}\text{ be} \\
h &=&\mathrm{dist}_{p}\circ \pi \circ \exp _{\tilde{p}}
\end{eqnarray*}%
From the proof of Proposition \ref{linear inf harm morph}, $f$ is not
infinity-harmonic. The $C^{1}$--distance between $f|_{B\left( 0,r\right) }$
and $h|_{B\left( 0,r\right) }$ goes to $0$ like $O\left( r^{2}\right) $ as $%
r\rightarrow 0$; so $h$ is not infinity-harmonic.

Combining the facts that $h$ is not infinity-harmonic, $\exp _{\tilde{p}%
}^{-1}$ is an $O\left( r^{2}\right) $--isometry on $B\left( 0,r\right) ,$
and 
\begin{equation*}
\left( \exp _{\tilde{p}}^{-1}\right) ^{\ast }\left( h\right) =\pi ^{\ast
}\left( \mathrm{dist}\left( p,\cdot \right) \right)
\end{equation*}%
we see that $\pi ^{\ast }\left( \mathrm{dist}_{p}\right) $ is not infinity-harmonic as desired.

The reader may be concerned that we \textquotedblleft run out of
room\textquotedblright\ for this argument, since we have to take $r$ very
small to make it work. This is not a concern, since $\Delta _{\infty }\left(
f\right) $ becomes arbitrarily large (in places) near the origin.

Now suppose that $d\pi _{\tilde{p}}\ $is nonzero and not onto. The above
argument shows that it is horizontally weakly conformal. As in the proof of
Proposition \ref{linear inf harm morph}, we take $v\in T_{p}B,$ to be
perpendicular to $\mathrm{Im}\left[ d\pi _{\tilde{p}}\right] $ and very
small. We saw that the $\infty$-laplacians of $\mathrm{dist}\left( v,\cdot
\right) $ and 
\begin{eqnarray*}
f :T_{\tilde{p}}E\longrightarrow \mathbb{R}, \;\;\; f =\mathrm{dist}%
_{v}\circ d\pi _{\tilde{p}}
\end{eqnarray*}%
both can be made arbitrarily large by choosing the norm of $v$ to be small
enough. Now let 
\begin{eqnarray*}
h :T_{\tilde{p}}E\longrightarrow \mathbb{R}\text{ be} \;\;\; h =\mathrm{dist}%
_{v}\circ \exp _{p}^{-1}\circ \pi \circ \exp _{\tilde{p}}.
\end{eqnarray*}%
The $C^{1}$--distance between $f|_{B\left( 0,r\right) }$ and $h|_{B\left(
0,r\right) }$ is $O\left( r^{2}\right) $ as $r\rightarrow 0$; so $h$ is not
infinity-harmonic. Since 
\begin{equation*}
\left( \exp _{\tilde{p}}^{-1}\right) ^{\ast }\left( h\right) =\pi ^{\ast
}\left( \mathrm{dist}_{v}\circ \exp _{p}^{-1}\right) ,
\end{equation*}%
and $\exp _{\tilde{p}}^{-1}$ is an $O\left( r^{2}\right) $--isometry on $%
B\left( 0,r\right) $ we see that $\pi ^{\ast }\left( \mathrm{dist}_{v}\circ
\exp _{p}^{-1}\right) $ is not infinity-harmonic, and in fact has $\infty$-laplacians that are arbitrarily large if the norm of $v$ is small enough.

Finally notice that the $C^{1}$--distance between $\pi ^{\ast }\left( 
\mathrm{dist}_{v}\circ \exp _{p}^{-1}\right) $ and $\pi ^{\ast }\left(
\left( \mathrm{dist}_{\exp _{p}\left( v\right) }\right) \right) $ converges
to $0$ as the norm of $v$ goes to zero. So $\pi ^{\ast }\left( \left( 
\mathrm{dist}_{\exp _{p}\left( v\right) }\right) \right) $ can not be
infinity-harmonic, even though $\mathrm{dist}_{\exp _{p}\left( v\right) }$
is infinity-harmonic. So infinity-harmonic morphisms are horizontally weakly
conformal maps.
\end{proof}

\begin{lemma}
\label{HWC2}If $\pi :E\rightarrow B$ is an infinity-harmonic morphism, then $%
\pi $ is an infinity-harmonic map.
\end{lemma}

\begin{proof}
By Lemma \ref{toHWC}, $\pi $ is a horizontally weakly conformal map. Let $%
\lambda $ be the dilation of $\pi $. Then, for any function $f$ locally
defined on $B$, we have 
\begin{eqnarray}
\left\vert \nabla (f\circ \pi )\right\vert ^{2} &=&g^{ij}(f\circ \pi
)_{i}(f\circ \pi )_{j}=g^{ij}f_{\alpha }\pi _{i}^{\alpha }f_{\beta }\pi
_{j}^{\beta }  \notag  \label{LDA} \\
&=&\lambda ^{2}(h^{{\alpha \beta }}\circ \pi )f_{\alpha }f_{\beta }=\lambda
^{2}(\left\vert \nabla f\right\vert ^{2}\circ \pi ),  \notag
\end{eqnarray}%
where the third equality was obtained by using the horizontally weakly
conformal equation $g^{ij}\pi _{i}^{\alpha }\pi _{j}^{\beta }=\lambda
^{2}(h^{{\alpha \beta }}\circ \pi )$ . It follows that 
\begin{equation}
\nabla (\left\vert \nabla (f\circ \pi )\right\vert ^{2})=(\nabla \,\lambda
^{2})(\left\vert \nabla \,f\right\vert ^{2}\circ \pi )+\lambda ^{2}\nabla
\,(\left\vert \nabla \,f\right\vert ^{2}\circ \pi ),  \notag  \label{E17}
\end{equation}%
and hence \addtocounter{theorem}{1} 
\begin{eqnarray}
\Delta _{\infty }^{M}(f\circ \pi ) &=&\frac{1}{2}\,\langle \nabla \,(f\circ
\pi ),\;\nabla \,\left\vert \nabla (f\circ \pi )\right\vert ^{2}\rangle 
\notag  \label{E30} \\
&=&\frac{1}{2}\,\langle \nabla \,(f\circ \pi ),\;(\nabla \,\lambda
^{2})(\left\vert \nabla \,f\right\vert ^{2}\circ \pi )\rangle
\label{Inf
Lapl f circ pi} \\
&+&\frac{1}{2}\,\langle \nabla \,(f\circ \pi ),\;\lambda ^{2}\nabla
\,(\left\vert \nabla \,f\right\vert ^{2}\circ \pi \rangle .  \notag \\
&=&\frac{1}{2}(\left\vert \nabla \,f\right\vert ^{2}\circ \pi )\mathrm{d}%
f\left( \mathrm{d}\pi (\nabla \lambda ^{2})\right)  \notag \\
&+&\frac{1}{2}\,\lambda ^{4}\langle \nabla \,f,\;\nabla \,\left\vert \nabla
f\right\vert ^{2}\rangle _{h}\circ \pi  \notag
\end{eqnarray}%
for any function (locally) defined on $B$.

Now choose $f$ to be a (locally defined) distance function. The second term
vanishes by the infinity-harmonicity of $f.$ Since $f\circ \pi $ is infinity-harmonic, 
\begin{eqnarray*}
0 &=&\Delta _{\infty }^{M}(f\circ \pi ) \\
&=&\frac{1}{2}\mathrm{d}f\left( \mathrm{d}\pi (\nabla \lambda ^{2})\right) ,
\end{eqnarray*}%
for \emph{any }locally defined distance function $f$ on $B.$ Therefore 
\begin{equation*}
\mathrm{d}\pi (\nabla \lambda ^{2})=0,
\end{equation*}%
and $\pi $ is a horizontally homothetic map. Applying Proposition \ref{P1}
we obtain the lemma.
\end{proof}

\begin{lemma}
\label{fromHWC} A horizontally weakly conformal, infinity-harmonic map is an
infinity-harmonic morphism.
\end{lemma}

\begin{proof}
Suppose $\pi :E\longrightarrow B$ is a horizontally weakly conformal
infinity-harmonic map, and $f:B\longrightarrow \mathbb{R}$ is any (locally
defined) function on $B$. It follows from Proposition \ref{P1} that $\pi $
is horizontally weakly conformal with dilation $\lambda $ having vertical
gradient. At points where $\pi $ is submersive, we have as in the proof of
Lemma \ref{HWC2} 
\begin{equation*}
\left\vert \mathrm{d}(f\circ \pi )\right\vert ^{2}=\lambda ^{2}(\left\vert 
\mathrm{d}\,f\right\vert ^{2}\circ \pi ).
\end{equation*}

At points where $\pi $ is critical, $d\pi $ is $0$ so $\left\vert \mathrm{d}%
(f \circ \pi )\right\vert ^{2}=0$. Since $\lambda $ is also zero at these
points we have 
\begin{equation*}
\left\vert \mathrm{d}(f \circ \pi )\right\vert ^{2}=\lambda ^{2}(\left\vert 
\mathrm{d}\,f \right\vert ^{2}\circ \pi ).
\end{equation*}%
in all cases.

Using Equation (\ref{Inf Lapl f circ pi}) and the fact that $\lambda $ has
vertical gradient we have 
\begin{equation*}
\Delta _{\infty }^{M}(f\circ \pi )=\frac{1}{2}\,\lambda ^{4}\langle \nabla
\,f,\;\nabla \,\left\vert \nabla f\right\vert ^{2}\rangle _{h}\circ \pi ,
\end{equation*}%
for any $f$ (locally) defined on $B$. This implies that $\pi $ pulls back
infinity-harmonic functions to infinity-harmonic functions and hence, by
definition, $\pi $ is an infinity-harmonic morphism.
\end{proof}

\begin{proof}
\emph{(of Theorem \ref{HWCTHM copy(1)}):} Together Lemmas \ref{toHWC}, \ref%
{HWC2}, \ref{fromHWC}, and Proposition \ref{inf harm subm} give Theorem \ref%
{HWCTHM copy(1)}.
\end{proof}

We conclude this section by pointing out that a proof similar to that of
Lemmas \ref{HWC2} and \ref{fromHWC} gives the following

\begin{proposition}
A map between Riemannian manifolds is an infinity-harmonic morphism if and
only if it pulls back locally defined infinity-harmonic maps to
infinity-harmonic maps.
\end{proposition}

\section{Constructions of Infinity-harmonic Maps}

In this section we give several methods to construct infinity-harmonic maps
into Euclidean space. We characterize those immersions which are
infinity-harmonic maps. We also show that isometrically immersing the target
manifold of a map into another manifold does not change the
infinity-harmonicity of the map. Coupled with Nash's embedding theorem this
suggests that it is particularly interesting to study the infinity-harmonic
maps into a Euclidean space.

\begin{proposition}
\label{R}$[Infinity-harmonic\,maps\,into\,a\,Euclidean\,space]$ \newline
A map $\varphi :\Omega \subset R^{m}\longrightarrow R^{n}$ with $\varphi
(x_{1},\ldots ,x_{m})=(\varphi ^{1}(x),\ldots ,\varphi ^{n}(x))$ is an
infinity-harmonic map if and only if it is a solution of the following
system of PDEs:%
\begin{equation*}
\langle \nabla \varphi ^{i},\nabla \left\vert \nabla \varphi ^{1}\right\vert
^{2}\rangle +\langle \nabla \varphi ^{i},\nabla \left\vert \nabla \varphi
^{2}\right\vert ^{2}\rangle +\cdots +\langle \nabla \varphi ^{i},\nabla
\left\vert \nabla \varphi ^{n}\right\vert ^{2}\rangle =0,\text{ for all }%
i=1,\ldots ,n
\end{equation*}
\end{proposition}

\begin{proof}
This is an easy exercise that we leave to the reader.
\end{proof}

\begin{proposition}
\label{ISO} Let $\varphi :(M,g)\longrightarrow (N,h)$ be a map and $\iota
:(N,h)\longrightarrow (Q,k)$ an isometric immersion. Then, the composition $%
\iota \circ \varphi :(M,g)\longrightarrow (Q,k)$ is infinity-harmonic if and
only if $\varphi $ is infinity-harmonic.
\end{proposition}

\begin{proof}
Since $\iota $ is an isometric immersion, 
\begin{equation}
\left\vert \mathrm{d}(\iota \circ \varphi )\right\vert ^{2}=\left\vert 
\mathrm{d}\varphi \right\vert ^{2}.  \notag
\end{equation}

So $\iota \circ \varphi $ is infinity-harmonic if and only if $\varphi $ is
infinity-harmonic.
\end{proof}

\begin{example}
Let $\varphi :R^{2}\setminus \{0\}\longrightarrow R^{2}$ be given by $%
\varphi (x,y)=(ax+by+c,\sqrt{x^{2}+y^{2}})$, where $a,b,c$ are constant.
Then $\varphi $ is an infinity-harmonic map with constant energy density.
Note that $\varphi $ is not an affine map, neither is it an isometric
immersion nor a Riemannian submersion.
\end{example}

It is well known that a map into Euclidean space is harmonic if and only if
its component functions are harmonic functions. So any choice of a set of
harmonic functions as components produces a harmonic map into Euclidean
space. It is easily checked that this is not true for infinity-harmonic maps
in general. Nevertheless, we have the following method to construct infinity
harmonic maps into Euclidean space using infinity-harmonic functions.

\begin{proposition}
\label{GZ4} Let $u:(M,g)\longrightarrow R$ and $v:(N,h)\longrightarrow R$ be
two infinity-harmonic functions. Then, $\varphi :(M\times
N,g+h)\longrightarrow R^{2}$ with $\varphi (x,y)=(u(x),v(y))$ for any $%
(x,y)\in M\times N$ is an infinity-harmonic map.
\end{proposition}

\begin{proof}
We leave the proof as an exercise.
\end{proof}

\begin{example}
It is well known \cite{Ar2} that $u(x_{1},x_{2})=x_{1}^{4/3}-x_{2}^{4/3}$ is
an infinity-harmonic function on $R^{2}$. By Proposition \ref{GZ4}, we have
a globally defined infinity-harmonic map $\varphi :R^{4}\longrightarrow
R^{2} $ given by $\varphi
(x_{1},x_{2},x_{3},x_{4})=(x_{1}^{4/3}-x_{2}^{4/3},%
\;x_{3}^{4/3}-x_{4}^{4/3}) $ which has nonconstant energy density $%
\left\vert \mathrm{d}\varphi \right\vert ^{2}=\frac{16}{9}%
(x_{1}^{2/3}+x_{2}^{2/3}+x_{3}^{2/3}+x_{4}^{2/3})$.
\end{example}

\begin{proposition}
$[Direct\,sum\,construction]$ Let $\varphi :(M,g)\longrightarrow R^{n}$ and $%
\psi :(N,h)\longrightarrow R^{n}$ be two infinity-harmonic maps into
Euclidean space. Then, their direct sum $\varphi \oplus \psi :(M\times
N,g+h)\longrightarrow R^{n}$ defined by $(\varphi \oplus \psi )(p,q)=\varphi
(p)+\psi (q)$ is an infinity-harmonic map.
\end{proposition}

\begin{proof}
We leave the proof as an exercise.
\end{proof}

\begin{proposition}
\label{00hm} An immersion $\iota :(M^{m},g)\longrightarrow (N^{n},h)$ is
infinity-harmonic if and only if the energy density of $\iota $ is constant.
In particular, a conformal immersion is infinity-harmonic if and only if it
is a homothetic immersion.
\end{proposition}

\begin{proof}
We leave the proof as an exercise.
\end{proof}

\begin{corollary}
\label{iden} Let $g$ and $h$ be two Riemannian metrics on a manifold $M$.
Then, the identity map $1:(M,g)\longrightarrow (M,h)$ is infinity-harmonic
if and only if $Trace_{g}h=constant$. In particular, $1:(M,g)\longrightarrow
(M,\lambda ^{2}g)$ is infinity-harmonic if and only if $\lambda $ is
constant.
\end{corollary}

\section{Infinity-harmonic Maps Into Spheres}

It is well known that in the presence of sufficient symmetry the harmonic
map equation can be reduced to an ordinary differential equation. In this
section we use ideas similar to those of Schoen and Yau's and Smith's about
harmonic maps to find infinity-harmonic maps into spheres.\newline

Let $i:R^{n}\longrightarrow R^{n+1}$ be the injection $i(x^{1},\ldots
,x^{n})=(x^{1},\ldots ,x^{n},0)$, and let $(r,\theta )$ denote the polar
coordinates on the unit ball $B^{n}$ and $(\rho ,\phi )$ the geodesic
coordinates on the unit sphere $S^{n}$, where $\rho $ is the distance from
the north pole of $S^{n}$ and $\phi \in S^{n-1}$. Following Schoen and Yau's
idea (\cite{SY}) we try to solve the Dirichlet Problem for rotationally
symmetric infinity-harmonic maps.

\begin{theorem}
A rotationally symmetric map $\varphi :B^{n}\longrightarrow S^{n}\subset
R^{n+1}$ of the form \addtocounter{theorem}{1} 
\begin{eqnarray}
\varphi &:&B^{n}\longrightarrow S^{n},  \label{symm} \\
\varphi (r,\theta ) &=&(\rho (r),\theta )\;\;with\;\;\rho (1)=\pi /2.  \notag
\end{eqnarray}
is an infinity-harmonic map if and only if either $\rho =\pi /2$ and $%
\varphi $ is the equator map 
\begin{eqnarray*}
\varphi &:&B^{n}\setminus \{0\}\longrightarrow S^{n-1}\subset S^{n}, \\
\varphi (x) &=&x/\left\vert x\right\vert ,
\end{eqnarray*}
or $\rho $ satisfies the ordinary differential equation %
\addtocounter{theorem}{1} 
\begin{equation}
{\rho ^{^{\prime }}}^{2}+\frac{n-1}{r^{2}}\sin ^{2}\rho =\mathrm{constant},
\label{O1}
\end{equation}%
and $\varphi $ has constant energy density.
\end{theorem}

\begin{proof}
Using polar coordinates $(r,\theta )$ on $B^{n}$ and geodesic coordinates $%
(\rho ,\phi )$ on $S^{n}$ we can write the metrics as 
\begin{equation}
g_{B^{n}}=\mathrm{d}r^{2}+r^{2}\mathrm{d}\theta ^{2},\;\;\;\mathrm{and}%
\;\;g_{S^{n}}=\mathrm{d}\rho ^{2}+\sin ^{2}\rho \,\mathrm{d}\phi ^{2}. 
\notag
\end{equation}%
Let $k_{ab}$ and $k^{ab}$ be the covariant and contravariant components of
the standard metric on $S^{n-1}$. Then, 
\begin{equation}
(g_{B^{n}})^{ij}=\left( 
\begin{array}{cc}
1 & 0 \\ 
0 & \frac{1}{r^{2}}k^{ab}%
\end{array}%
\right) ,\;\;(g_{S^{n}})_{\alpha \beta }\circ \varphi =\left( 
\begin{array}{cc}
1 & 0 \\ 
0 & \sin ^{2}\rho k_{ab}%
\end{array}%
\right) .  \notag
\end{equation}%
Since 
\begin{equation}
\mathrm{d}\varphi =(\varphi _{i}^{\alpha })=\left( 
\begin{array}{cc}
\rho ^{^{\prime }} & 0 \\ 
0 & id%
\end{array}%
\right) ,  \notag
\end{equation}%
we have 
\begin{equation}
\left\vert \mathrm{d}\varphi \right\vert ^{2}=(g_{B^{n}})^{ij}\varphi
_{i}^{\alpha }\varphi _{j}^{\beta }(g_{S^{n}})_{\alpha \beta }\circ \varphi =%
{\rho ^{^{\prime }}}^{2}+\frac{n-1}{r^{2}}\sin ^{2}\rho .  \notag
\end{equation}%
and 
\begin{equation*}
\nabla \left\vert \mathrm{d}\varphi \right\vert ^{2}=\left( {\rho ^{^{\prime
}}}^{2}+\frac{n-1}{r^{2}}\sin ^{2}\rho \right) ^{\prime }\partial _{r}
\end{equation*}

So $\varphi $ is infinity-harmonic if and only if 
\begin{eqnarray*}
0 &=&d\varphi \left( \nabla \left\vert \mathrm{d}\varphi \right\vert
^{2}\right) \\
&=&\rho ^{^{\prime }}\left( {\rho ^{^{\prime }}}^{2}+\frac{n-1}{r^{2}}\sin
^{2}\rho \right) ^{\prime }.
\end{eqnarray*}

It follows that either $\rho =constant$ and hence $\rho =\pi /2$ by boundary
condition, or $\rho $ is a solution of the ODE 
\begin{equation}
{\rho ^{^{\prime }}}^{2}+\frac{n-1}{r^{2}}\sin ^{2}\rho =\mathrm{constant}. 
\notag
\end{equation}%
The first case corresponds to the map $\varphi (r,\theta )=(\pi /2,\theta )$
which, in Cartesian coordinates, can be expressed as $\varphi
:B^{n}\setminus \{0\}\longrightarrow S^{n-1}\subset S^{n},\;\;\varphi
(x)=x/\left\vert x\right\vert $, the equator map. In the second case we have 
\begin{equation}
\left\vert d\varphi \right\vert ^{2}={\rho ^{^{\prime }}}^{2}+\frac{n-1}{%
r^{2}}\sin ^{2}\rho =\mathrm{constant}  \notag
\end{equation}%
as desired.
\end{proof}

Let $S^{2}$ be the unit sphere in $R^{3}$ parametrized by spherical polar
coordinates: 
\begin{equation}
(\alpha ,\beta )\longmapsto (\cos \alpha ,\sin \alpha e^{i\beta })\in 
\mathbb{R}\oplus \mathbb{C},\text{ }(\alpha ,\beta )\in \mathbb{R}\times 
\mathbb{R}  \notag
\end{equation}%
Parametrized the cylinder as $\mathbb{R}\times S^{1}=\left\{
(s,e^{it}):(s,t)\in \mathbb{R}\times \mathbb{R}\right\} $, and consider the
rotationally symmetric map $\varphi :R\times S^{1}\longrightarrow S^{2}$ of
the form \addtocounter{theorem}{1} 
\begin{equation}
\varphi (s,t)=(\cos \alpha (s),\sin \alpha (s)e^{ikt}),  \label{TORUS}
\end{equation}%
where $\alpha :R\longrightarrow R$ is a smooth function and $k$ a non-zero
integer. Smith (\cite{Sm}) proved that $\varphi $ is harmonic if and only if 
$\alpha $ is a solution of the ordinary differential equation %
\addtocounter{theorem}{1} 
\begin{equation}
\alpha ^{\prime \prime }=k^{2}\sin \alpha \cos \alpha ,  \label{SMITH}
\end{equation}%
and by solving this equation of pendulum with constant gravity and no
damping he was able to find some interesting harmonic maps from torus into a
sphere (see also \cite{BW1} for a detailed explanation).\newline

Our next theorem shows that Smith's method can also be used to find \newline
infinity-harmonic maps into a sphere.

\begin{theorem}
The rotationally symmetric map $\varphi :R\times S^{1}\longrightarrow S^{2}$%
, 
\begin{equation}
\varphi (s,t)=(\cos \alpha (s),\sin \alpha (s)e^{ikt})  \notag
\end{equation}%
is an infinity-harmonic map if and only if

\begin{itemize}
\item[(1)] $\alpha =$ constant and $\varphi $ is the projection onto the
second factor followed by a homothetic immersion, or

\item[(2)] $\alpha (s)=2\arctan (e^{ks+A})-\pi /2$, where $A$ is any
constant, or

\item[(3)] $2\alpha $ is a solution of the pendulum equation (Equation 28.74
in \cite{TP}) \addtocounter{theorem}{1} 
\begin{equation}
\frac{\mathrm{d}^{2}\theta }{\mathrm{d\,t}^{2}}+k^{2}\sin \theta =0.
\label{PDULUM}
\end{equation}%
In this case, the map $\varphi $ factors to an infinity-harmonic map from
the torus\newline
$R/\langle T\rangle \times S^{1}$ to $S^{2}$, where $T$ is the period of $%
\alpha $.
\end{itemize}
\end{theorem}

Note that the ODE (\ref{PDULUM}) differs from the harmonic map equation (\ref%
{SMITH}) by a negative sign.

\begin{proof}
The infinity-harmonicity of $\varphi :R\times S^{1}\longrightarrow S^{2}$ is
the same as the infinity-harmonicity of $\varphi :R\times
S^{1}\longrightarrow S^{2}\hookrightarrow R^{3}$. Changing from complex to
real notation we have 
\begin{equation*}
\varphi (s,t)=(\cos \alpha (s),\sin \alpha (s)e^{ikt})=(\cos \alpha
(s),\;\sin \alpha (s)\cos kt,\;\sin \alpha (s)\sin kt),
\end{equation*}%
so \addtocounter{theorem}{1} 
\begin{equation}
\begin{cases}
& \nabla \varphi ^{1}=(-\alpha ^{\prime }(s)\sin \alpha (s),\;0), \\ 
& \nabla \varphi ^{2}=(\alpha ^{\prime }(s)\cos \alpha (s)\cos kt,\;-k\sin
\alpha (s)\sin kt), \\ 
& \nabla \varphi ^{3}=(\alpha ^{\prime }(s)\cos \alpha (s)\sin kt,\;k\sin
\alpha (s)\cos kt),\;\mathrm{and} \\ 
& \left\vert \mathrm{d}\varphi \right\vert ^{2}=\sum_{i=1}^{3}\left\vert
\nabla \varphi ^{i}\right\vert ^{2}={\alpha ^{\prime }}(s)^{2}+k^{2}\sin
^{2}\alpha (s).%
\end{cases}%
\end{equation}%
By Proposition \ref{R}, $\varphi $ is infinity-harmonic if and only if%
\begin{equation*}
\begin{cases}
& -\alpha ^{\prime }(s)\sin \alpha (s)({\alpha ^{\prime }}^{2}+k^{2}\sin
^{2}\alpha )^{\prime }=0 \\ 
& \alpha ^{\prime }(s)\cos \alpha (s)\cos kt({\alpha ^{\prime }}%
^{2}+k^{2}\sin ^{2}\alpha )^{\prime }=0 \\ 
& \alpha ^{\prime }(s)\cos \alpha (s)\sin kt({\alpha ^{\prime }}%
^{2}+k^{2}\sin ^{2}\alpha )^{\prime }=0.%
\end{cases}%
\end{equation*}%
It follows that $\varphi $ is infinity-harmonic if and only if either $%
\alpha ^{\prime }=0$ and hence $\alpha $ is constant and $\varphi $ is the
projection onto the second factor followed by a homothetic immersion, or, 
\begin{equation*}
({\alpha ^{\prime }}^{2}+k^{2}\sin ^{2}\alpha )^{\prime }=0,
\end{equation*}%
which is equivalent to 
\begin{equation*}
{\alpha ^{\prime }}^{2}+k^{2}\sin ^{2}\alpha =C.
\end{equation*}%
When $C=k^{2}$, we can solve the previous equation and get\newline
$\alpha =2\arctan (e^{ks+A})-\pi /2$. When $C>k^{2}$, 
\begin{equation*}
{\alpha ^{\prime }}^{2}+k^{2}\sin ^{2}\alpha =C.
\end{equation*}%
is equivalent to 
\begin{equation*}
(2\alpha )^{\prime \prime }+k^{2}\sin (2\alpha )=0,
\end{equation*}%
which means that $2\alpha $ is a solution of the pendulum equation (\ref%
{PDULUM}). By the theory of the solutions of the pendulum equation (see,
e.g., \cite{TP}) we obtain the statement (3).
\end{proof}

\section{The Effect of a Conformal Change on the Infinity-Laplacian}

In this section we study the effect of a conformal change on the $\infty $%
-Laplacian to derive formulas for the $\infty $-Laplacian of spheres and
hyperbolic spaces in terms of the $\infty $-Laplacian on Euclidean space.

Given Riemannian metrics $g$ and $\bar{g}$ on a smooth manifold $M,$ we let $%
\nabla $, $\left\vert .\right\vert $, and $\Delta _{\infty }$ denote the
gradient, the norm, and the $\infty $-Laplacian with respect to $g$ and we
let $\bar{\nabla}$, $\left\vert .\right\vert _{\bar{g}}$, and $\bar{\Delta}%
_{\infty }$ denote the gradient, the norm, and the $\infty $-Laplacian with
respect to $\bar{g}$.

\begin{theorem}
\label{CONF} Let $\bar{g}=F^{-2}g$ be a metric conformal to $g$ on $M$. Then %
\addtocounter{theorem}{1} 
\begin{equation}
\bar{\Delta}_{\infty }\,u=F^{4}\Delta _{\infty }\,u+F^{3}\left\vert \nabla
\,u\right\vert ^{2}g(\nabla \,u,\nabla \,F).
\end{equation}
\end{theorem}

\begin{proof}
A direct computation using $\bar{g}=F^{-2}g$ gives 
\begin{eqnarray*}
&&\bar{\nabla}\,u=F^{2}\nabla \,u, \\
\bar{\nabla}\,\left\vert \bar{\nabla}\,u\right\vert _{\bar{g}} ^{2} &=&\bar{%
\nabla}\,\left\vert F^{2}\nabla \,u\right\vert _{\bar{g}}^{2} =\bar{\nabla}%
\,\left( F^{4}F^{-2}\left\vert \nabla \,u\right\vert _{g}^{2}\right) \\
&=&2F^{3}\left\vert \nabla \,u\right\vert _{g}^{2}\nabla F+F^{4}\nabla
\left\vert \nabla \,u\right\vert _{g}^{2}
\end{eqnarray*}%
It follows that \addtocounter{theorem}{1} 
\begin{eqnarray}
\bar{\Delta}_{\infty }\,u &=&\frac{1}{2}\bar{g}(\bar{\nabla}\,u,\,\bar{\nabla%
}\,\left\vert \bar{\nabla}\,u\right\vert _{\bar{g}} ^{2})  \notag \\
&=&\frac{1}{2}F^{-2}g(F^{2}\nabla \,u,\,2F^{3}\left\vert \nabla
\,u\right\vert _{g}^{2}\nabla F+F^{4}\nabla \left\vert \nabla \,u\right\vert
_{g}^{2}) \\
&=&F^{4}\Delta _{\infty }\,u+F^{3}\left\vert \nabla \,u\right\vert
^{2}g(\nabla \,u,\nabla \,F).  \notag
\end{eqnarray}
\end{proof}

As an application of Theorem \ref{CONF} we have the following expression for
the $\infty$-Laplace equation in hyperbolic\thinspace space.

\begin{corollary}
\label{CHYPER} $[\infty -Laplacian\,on\,hyperbolic\,space\,B^{m}]$ Let $%
(B^{m},g^{H})$ be the $m$-dimensional hyperbolic space with open-ball model,
where $B^{m}=\{x\in R^{m}:\left\vert x\right\vert <1\}$ and $%
g^{H}=F^{-2}\delta _{ij}$ with $F=\frac{1}{2}(1-\left\vert x\right\vert
^{2}) $. Then, the $\infty $-Laplace equation in the hyperbolic space $%
(B^{m},g^{H})$ is the conformal $\infty $-Laplace equation in the Euclidean
space $(R^{m},\delta _{ij})$, which can be written as %
\addtocounter{theorem}{1} 
\begin{equation}
\Delta _{\infty }^{\mathbb{R}^{m}}\,u-\frac{2\left\vert \nabla
\,u\right\vert ^{2}}{1-\left\vert x\right\vert ^{2}}\langle x,\nabla
u\rangle =0,\;\;x\in \mathbb{R}^{m},  \label{hyper1}
\end{equation}%
where $|,|$ and $\nabla $ denote the norm and the gradient defined by the
Euclidean metric $\langle ,\rangle $ on $B^{m}\subset R^{m}$.
\end{corollary}

\begin{example}
Let $u:\Omega \subset (S^{2}\setminus \{N\},g_{can})\cong
(R^{2},F^{-2}\delta _{ij})\longrightarrow R$ be given by $%
u(x_{1},x_{2})=\arctan \frac{x_{1}}{x_{2}}$. Then, we know (see \cite{Ar1})
that $u$ is an infinity-harmonic function on $\Omega \subset R^{2}$, so $%
\Delta _{\infty }^{\mathbb{R}^{2}}\,u=0$. On the other hand, we can easily
check that $\langle x,\nabla u\rangle =0$. Therefore, $u$ satisfies Equation
(\ref{hyper1}) and hence it is an infinity-harmonic function on sphere $%
S^{2} $. A more geometric way to see this is via the isometric $\mathbb{R}$%
--action that rotates the $2$--sphere and Proposition \ref{isometric grp
actions}.

Note that the function $u(x_{1},x_{2})=\arctan \frac{x_{1}}{x_{2}}$ is also
an infinity-harmonic function on hyperbolic space $(B^{2},g^{H})$ wherever
it is defined.
\end{example}

The following example give families of infinity-harmonic functions on
hyperbolic space.

\begin{example}
Let $(B^{m},g^{H})$ be the $m$-dimensional hyperbolic space with open-ball
model as in Corollary \ref{CHYPER}. Then, for constants $a_{1},...,a_{m-1}$,
the function $u:(B^{m},g^{H})\longrightarrow R$ given by $%
u(x)=(a_{1}x_{1}+...+a_{m-1}x_{m-1})(1+\left\vert x\right\vert
^{2}-2x_{m})^{-1}$ is an infinity-harmonic function. This follows from
Theorem 2.9 in \cite{Ou}.
\end{example}

\end{document}